\documentclass[a4paper,12pt]{amsart}

\usepackage{geometry}
\usepackage{lmodern}
\usepackage[stretch=10]{microtype}
\usepackage{graphicx}
\usepackage{tikz-cd}

\usepackage{amsmath}
\usepackage{amssymb}
\usepackage{amsthm}
\usepackage[matrix,arrow,cmtip]{xy}
\usepackage{bbm}
\usepackage[shortlabels]{enumitem}
\PassOptionsToPackage{hyphens}{url}
\usepackage[bookmarksnumbered,bookmarksdepth=2,bookmarksopen,colorlinks,linkcolor=black,citecolor=black,urlcolor=black,linktoc=all]{hyperref}
\usepackage[capitalise]{cleveref}

\newcommand{\bA}{\mathbb{A}}
\newcommand{\bC}{\mathbb{C}}
\newcommand{\bF}{\mathbb{F}}
\newcommand{\bG}{\mathbb{G}}
\newcommand{\bH}{\mathbb{H}}

\newcommand{\bN}{\mathbb{N}}
\newcommand{\bP}{\mathbb{P}}
\newcommand{\bQ}{\mathbb{Q}}
\newcommand{\bR}{\mathbb{R}}
\newcommand{\bS}{\mathbb{S}}

\newcommand{\bZ}{\mathbb{Z}}

\newcommand{\cA}{\mathcal{A}}

\newcommand{\cC}{\mathcal{C}}
\newcommand{\cD}{\mathcal{D}}

\newcommand{\cF}{\mathcal{F}}
\newcommand{\cH}{\mathcal{H}}

\newcommand{\cM}{\mathcal{M}}
\newcommand{\cO}{\mathcal{O}}

\newcommand{\cS}{\mathcal{S}}
\newcommand{\cT}{\mathcal{T}}
\newcommand{\cV}{\mathcal{V}}
\newcommand{\cX}{\mathcal{X}}

\newcommand{\gG}{\mathbf{G}}
\newcommand{\gGL}{\mathbf{GL}}
\newcommand{\gGSp}{\mathbf{GSp}}

\newcommand{\gH}{\mathbf{H}}

\newcommand{\gM}{\mathbf{M}}
\newcommand{\gMT}{\mathbf{MT}}

\newcommand{\gP}{\mathbf{P}}

\newcommand{\gQ}{\mathbf{Q}}

\newcommand{\gSL}{\mathbf{SL}}

\newcommand{\gSp}{\mathbf{Sp}}

\newcommand{\gU}{\mathbf{U}}
\newcommand{\gX}{\mathbf{X}}

\newcommand{\gZ}{\mathbf{Z}}

\newcommand{\rM}{\mathrm{M}}

\DeclareMathOperator{\Aut}{Aut}
\DeclareMathOperator{\Atyp}{Atyp}

\DeclareMathOperator{\disc}{disc}

\DeclareMathOperator{\Gal}{Gal}

\DeclareMathOperator{\Opt}{Opt}
\DeclareMathOperator{\Res}{Res}

\DeclareMathOperator{\Sh}{Sh}

\DeclareMathOperator{\Stab}{Stab}

\newcommand{\ad}{\mathrm{ad}}

\newcommand{\der}{\mathrm{der}}

\newcommand{\fullsmallmatrix}[4]{\bigl( \begin{smallmatrix} #1 & #2 \\ #3 & #4 \end{smallmatrix} \bigr)}

\newtheorem{lemma}{Lemma}[subsection]

\newtheorem{theorem}[lemma]{Theorem}

\newtheorem{conjecture}[lemma]{Conjecture}

\Crefname{conjecture}{Conjecture}{Conjectures} 

\Crefname{claim}{Claim}{Claims}

\newtheorem{remark}[lemma]{Remark}
\newtheorem{example}[lemma]{Example}

\newtheorem*{lemma*}{Lemma}
\newtheorem*{proposition*}{Proposition}
\newtheorem*{theorem*}{Theorem}
\newtheorem*{corollary*}{Corollary}
\newtheorem*{claim*}{Claim}

\theoremstyle{definition}
\newtheorem*{definition}{Definition}

\usepackage{ifthen}
\newcounter{constant}

\newcommand{\newC}[1]{%
   \ifthenelse{\equal{#1}{*}} {%
      \stepcounter{constant} c_{\theconstant}%
   } {%
      \refstepcounter{constant} c_{\theconstant} \label{C:#1}%
   }%
}
\newcommand{\refC}[1]{c_{\ref*{C:#1}}}

\usepackage{color}

\newcommand{\Qbar}{\overline{\bQ}}

\usepackage{csquotes}

\title[Unlikely intersections in Shimura varieties]{Unlikely intersections in Shimura varieties and beyond: a survey}

\author{Christopher Daw}

\address{}
\email{}

\dedicatory{}

\subjclass[2020]{}

\begin{document}

\begin{abstract}
The aim of this note is to provide a concise introduction to so-called problems of unlikely intersections for (pure) Shimura varieties and to review the current state of the art. In the process, we will touch upon more general settings and some of the results in those contexts.
\end{abstract}

\maketitle

\tableofcontents

\section{Introduction}

The term {\it unlikely intersections} seems to have originated in the collaborations of Bombieri--Masser--Zannier around the turn of the millennium (see \cite[Th. 2]{BMZ99} for an early incarnation and \cite{BMZ08} for an early article making explicit reference to the term). It is also closely related to other terminologies, such as Zilber's {\it atypical} intersections \cite[Def. 2]{Zilber02} and Bombieri--Masser--Zannier's {\it anomalous subvarieties} \cite[Def. 1.1]{BMZ:anom}.

Through the works of Pink \cite{Pin05,pink:generalisation} and others, it has come to unify a wealth of results and conjectures across arithmetic geometry, some of which are far older than it---prototypical in this regard is the Manin--Mumford conjecture, for instance. In essence, it is a principle, which can be described as follows.

\begin{displayquote}
{\it
    Let $X$ be an algebraic variety endowed with a countable collection of distinguished algebraic subvarieties. Pick an arbitrary irreducible subvariety $V$ of $X$ and consider the intersection of $V$ with a distinguished subvariety of codimension exceeding $\dim V$. If $V$ is picked at random, this intersection is most likely empty. Therefore, the (countable) union of all such intersections ought to be sparse in $V$. 
    }
\end{displayquote}

Of course, if $V$ is not generic, the principle is faulty---if $V$ is contained in a proper subvariety of $X$, which also contains an abundance of distinguished subvarieties, there may well be an abundance of intersections. The field of unlikely intersections is, in large part, an attempt to prove, in an array of different settings (algebraic tori, abelian varieties, Shimura varieties, variations of Hodge structures...), that this is the only way in which the principle can fail (in characteristic $0$, at least). 

There are many brilliant surveys of this topic. The book of Zannier \cite{Zan12} has played a key role in popularising the field. The collection \cite{HRSUY:CIRM}, appearing at a similar time, contains several valuable expositions. Pila's book \cite{pila:book} is wide-ranging and contains many of the most recent advances. The focus of this article, therefore, will be specifically on unlikely intersections in pure Shimura varieties and the current state-of-the-art in this setting. There will be a particular emphasis on the arithmetic aspects, where there is still a vast amount to do. There are two primary intended audiences:
\begin{itemize}
    \item[(1)] early career researchers (particularly PhD students) who aspire to work in the area;
    \item[(2)] experts in other (most likely, nearby and related) fields, who are looking for an efficient summary of the concepts, techniques, results and conjectures.
\end{itemize}

In \cite{Daw15}, we wrote an article not entirely dissimilar in style to the present one, focusing on the so-called {\it Pila--Zannier strategy} for the Andr\'e--Oort conjecture. Klingler--Ullmo--Yafaev subsequently wrote a more contemporary survey of this material in \cite{KUY:bialg}. Below, there will certainly be some overlap with these texts. However, in \cite{Daw15}, we gave a much more detailed (and somewhat narrational) account of Shimura varieties, whereas, here, we will be much more concise. Similarly, the descriptions of Pila--Zannier were more detailed in \cite{Daw15} and \cite{KUY:bialg} but restricted to Andr\'e--Oort. Here we will be more brief, in order to cover more ground---a decade later, the field is considerably more advanced.

We will comment briefly on topics beyond our primary subject matter, such as unlikely intersections for mixed Shimura varieties and variations of Hodge structures. We will not go into detail, however, as there are also several excellent recent surveys in this direction (see, for example, \cite{Klingler:ICM}, \cite{Baldi:survey} and \cite{baldi:hab}). Hopefully, our omissions are not overly egregious and the references we include herein will prove sufficient to lead the interested reader in the right directions. 

Unfortunately, we will also say very little on many related and vitally important topics---such as the theory of o-minimal structures and the Pila--Wilkie theorem, methods from differential algebra, functional transcendence, equidistribution, effectivity, and results in characteristic $p$. Fortunately, these are covered in many other places (including the sources above and below).

\subsection*{Acknowledgements}

The author would like to express his deepest thanks to the organisers of the third Journal of Number Theory Biennial Conference held in Cetraro, Italy in August 2024. In addition, he wishes to thank Dorian Goldfeld and Federico Pellarin for managing the special volume associated with the conference. Finally, he heartily thanks Gregorio Baldi, Martin Orr, Jonathan Pila, Ziyang Gao, Umberto Zannier, Boris Zilber and the anonymous referee for many insightful comments that have surely helped to address some of the shortcomings in earlier drafts.

ChatGPT was used to proofread a final draft.

\section{Background material}

\subsection{Algebraic varieties}
 For an algebraic variety $V$ over a field $K$, and $L$ an extension of $K$, we denote by $V_L$ the base change of $V$ to $L$. If $K$ is a subfield of $\bC$, we denote by $V^{\rm an}$ the analytification of $V_\bC$. By an irreducible (sub)variety we refer to a geometrically irreducible (sub)variety, unless we say otherwise.

\subsection{Algebraic groups} 
For an algebraic group $\gG$, we denote by $\gG^\circ$ the (Zariski) connected component of $\gG$ containing the identity. Reductive and semisimple algebraic groups are by definition connected. We denote by $\gZ(\gG)$ the centre of $\gG$ and, if $\gH$ is an algebraic subgroup of $\gG$, we denote by $\gZ_{\gG}(\gH)$ the centraliser of $\gH$ in $\gG$.

For a reductive algebraic group $\gG$ we denote by $\gG^\der$ its derived (or commutator) subgroup and by $\gG^\ad$ the quotient $\gG/\gZ(\gG)$. Both $\gG^\der$ and $\gG^\ad$ are semisimple, and $\gG$ is the almost direct product $\gZ(\gG)\gG^{\rm der}$. When $\gG$ is defined over $\bR$, we denote by $\gG(\bR)^+$ the connected component of $\gG(\bR)$ containing the identity, and by $\gG(\bR)_+$ the preimage of $\gG^{\ad}(\bR)^+$ under the natural map $\gG(\bR)\to\gG^\ad(\bR)$. When $\gG$ is defined over $\bQ$, we write $\gG(\bQ)_+$ for $\gG(\bQ)\cap \gG(\bR)_+$.

\subsection{Arithmetic subgroups}
Let $\gG$ denote a reductive algebraic group over $\bQ$ and let $\gG\to\gGL_n$ be an embedding (also over $\bQ$). By an {\it arithmetic subgroup} of $\gG(\bQ)$, we refer to a subgroup that is commensurable with $\gG(\bQ)\cap\gGL_n(\bZ)$. The definition is independent of the embedding \cite[7.13]{Bor69}. Any arithmetic subgroup of $\gG(\bQ)$ is clearly a discrete subgroup of $\gG(\bR)$.

\subsection{Congruence subgroups}
By a {\it congruence subgroup} of $\gG(\bQ)$, we refer to a subgroup $\Gamma$ such that, for some $N\in\mathbb{N}$, the group
\[\Gamma(N):=\gG(\bQ)\cap\{\gamma\in\gGL_n(\bZ):\gamma\equiv{\rm id}\text{ mod }N\}\]
is a subgroup of $\Gamma$ of finite index.
Again, the definition is independent of the embedding. Indeed, if $K$ is a compact open subgroup of $\gG(\bA_f)$, then $K\cap\gG(\bQ)$ is a congruence subgroup of $\gG(\bQ)$ and every congruence subgroup of $\gG(\bQ)$ is of this form \cite[Prop. 4.1]{Mil05}. Clearly, a congruence subgroup of $\gG(\bQ)$ is an arithmetic subgroup of $\gG(\bQ)$.

\subsection{Hermitian symmetric domains}
By a {\it hermitian symmetric domain}, we refer to a connected, symmetric, hermitian manifold of non-compact type (see \cite[Ch. 1]{Mil05} for more details). For any hermitian symmetric domain $D$, we denote by ${\Aut}(D)$ the group of holomorphic isometries. The compact open topology on $D$ induces a canonical real Lie group structure on ${\Aut}(D)$. Moreover, for any $x\in D$, the stabiliser $K$ of $x$ in ${\Aut}(D)^+$ is maximal compact and ${\Aut}(D)^+/K\to D$ is an isomorphism of smooth manifolds (cf. \cite[Lem. 1.5]{Mil05}).

\subsection{Shimura data}
We write $\bG_m$ for the multiplicative group---that is, $\bG_m$ is the affine group scheme over $\bZ$ for which $\bG_m(R)=R^\times$ for any ring $R$---and we write $\bS$ for the {\it Deligne torus}, defined as the Weil restriction $\Res_{\bC/\bR}\bG_{m,\bC}$ (hence, by definition, $\bS(\bR)=\bC^\times$).

Given a reductive algebraic group $\gG$ over $\bQ$ and the $\gG(\bR)$-conjugacy class $\gX$ of a homomorphism $x_0:\bS\to\gG_{\bR}$ (where $(gx_0g^{-1})(s):=gx_0(s)g^{-1}$ for any $g\in\gG(\bR)$ and $s\in\bS(\bC)$), we say that the pair $(\gG,\gX)$ is a {\it Shimura datum} if it satisfies the properties SV1--SV3 of \cite[p. 302]{Mil05}. The primary significance of these properties is as follows.

\subsubsection{Relation with hermitian symmetric domains}  Let $X$ be a connected component of $\gX$. Then $X$ is naturally endowed with the structure of a hermitian symmetric domain for which the group $\gG(\bR)_+$ acts by holomorphic isometries, and the homomorphism
    \[\gG^\der(\bR)^+\to{\Aut}(X)^+\]
    is surjective with compact kernel \cite[Prop. 4.8 and Cor. 5.8]{Mil05}. Note that this homomorphism factors through $\gG^\der(\bR)^+\to\gG^\ad(\bR)^+$ (which is surjective by \cite[Prop. 5.1]{Mil05}). Moreover, any pair $(D,\gG_D)$ consisting of a hermitian symmetric domain $D$ and a semisimple algebraic group $\gG_D$ over $\bQ$ for which there exists a surjective homomorphism $\gG_D(\bR)^+\to\Aut(D)^+$ with compact kernel arises in this fashion.

\subsubsection{Realisations}

A {\it realisation} of $X$ is an analytic subset $\cX$ of a complex quasi-projective algebraic variety $\tilde\cX$ possessing a transitive holomorphic action of $\gG^{\rm der}(\bR)^+$ such that, for any $x_0\in\cX$, the map $\gG^{\rm der}(\bR)^+\to\cX$ sending $g$ to $gx_0$ is semi-algebraic and identifies $\cX$ with the quotient of $\gG^\der(\bR)^+$ by a maximal compact subgroup (and, therefore, with $X$). 

Two notable realisations arise from the {\it Harish--Chandra embedding}, which realises $X$ as a bounded symmetric domain in $\bC^n$ (for some $n\in\bN$), and the {\it Borel embedding}, which realises $X$ as a subset of the projective algebraic variety $X^\vee:=\gG(\bC)/\gP(\bC)$ with $\gP$ a parabolic subgroup of $\gG_\bC$. 

 \subsubsection{Algebraic structure}\label{sec:alg-structure}

The fact that $X$ has a bounded realisation shows that it is never an algebraic variety, unless it is a point. Nonetheless, there is a canonical algebraic structure on $X$. Indeed, let $\cX$ denote a realisation of $X$. A subset $Y\subset X$ is called an irreducible algebraic subvariety of $X$ if it corresponds to an irreducible analytic component of $\cX\cap \tilde{Y}$ for an algebraic subvariety $\tilde{Y}$ of $\tilde\cX$. This definition is independent of the choice of $\cX$ \cite[Cor. B.2]{kuy:ax-lindemann}.

Using the Borel embedding, we can endow $X$ with a $\Qbar$-algebraic structure. Indeed, given a faithful $\bQ$-algebraic representation $V$ of $\gG$, the variety $X^\vee$ is a subvariety of a (projective) flag variety $\Pi_\bC:=\gGL(V_\bC)/\gQ(\bC)$ with $\gQ$ a parabolic subgroup of $\gGL(V_\bC)$ containing $\gP$. The complex variety $\Pi_\bC$ has a natural $\bQ$-model $\Pi$ and, with respect to this model, the subvariety $X^\vee$ of $\Pi_\bC$ is defined over $\Qbar$ (see \cite[Sec. 3]{UY:characterisation} for more details).

\subsubsection{Decompositions} 

The group $\gG^{\rm ad}$ is equal to a direct product $\gG_1\times\cdots\times\gG_n$ of $\bQ$-simple subgroups. This induces a product decomposition $X=X_1\times\cdots\times X_n$ of hermitian symmetric domains. By grouping the factors, we obtain (finitely many) two-factor decompositions $\gG^{\rm ad}=\gG'_1\times\gG'_2$ and $X=X'_1\times X'_2$.

\subsubsection{Mumford--Tate groups}

For $x\in\gX$, there is a smallest (necessarily reductive) algebraic subgroup $\gH$ of $\gG$ defined over $\bQ$ such that $x(\bS)\subset\gH_\bR$. We refer to this group as the {\it Mumford--Tate group} of $x$, and we denote it $\gMT(x)$. It is easy to see that the set of $\gMT(x)$ for $x\in\gX$ contains a unique maximal element. We refer to this group as the Mumford--Tate group of $\gX$ and denote it $\gMT(\gX)$. Then $\gX$ is an $\gMT(\gX)(\bR)$-conjugacy class and, by \cite[Lem. 3.3]{ullmo:equidistribution}, $(\gMT(\gX),\gX)$ is a Shimura datum.

\subsection{Locally symmetric varieties}

Let $(\gG,\gX)$ be a Shimura datum and let $\Gamma$ be the intersection of an arithmetic subgroup of $\gG(\bQ)$ with $\gG(\bQ)_+$. By the theorem of Baily--Borel \cite[Th. 10.11]{BB:compactification}, for any connected component $X$ of $\gX$, the locally symmetric space $\Gamma\backslash X$ is a quasi-projective algebraic variety. Moreover, the algebraic structure on $\Gamma\backslash X$ is unique \cite[Cor. 3.16]{Mil05}. The embedding into projective space is given by automorphic forms on $X$ of sufficiently high weight, and we refer to the Zariski closure of the image as the {\it Baily--Borel compactification} of $\Gamma\backslash X$.

\subsubsection{Fundamental sets}

Let $K$ be a maximal compact subgroup of $\gG(\bR)$. By a {\it fundamental set} in $\gG(\bR)$ for $\Gamma$ and $K$, we refer to a subset $\cF_{\gG}$ of $\gG(\bR)$ such that $\cF_{\gG} K=\cF_{\gG}$, $\gG(\bR)=\Gamma\cF_{\gG}$ and, for any $g\in\gG(\bR)$, the set
\[\{\gamma\in\Gamma:\gamma\cF_{\gG}\cap g\cF_{\gG}\neq\emptyset\}\]
is finite (cf. \cite[Sec. 2D]{Orr18}). Fundamental sets were constructed using {\it Siegel sets} in \cite[Th. 13.1, 15.4]{Bor69}. We will henceforth restrict our definition to fundamental sets of this form. For the definition of a Siegel set, we appeal to \cite[Sec. 2.B]{Orr18}. 

If $x_0\in X$ satisfies $\Stab_{\gG(\bR)_+}(x_0)\subset K$, we refer to the image $\cF$ of $\cF_{\gG}\cap\gG(\bR)_+$ under the map $\gG(\bR)_+\to X$ given by $g\mapsto gx_0$ as
a fundamental set in $X$ for $\Gamma$. Clearly, $X=\Gamma\cF$ and the preimage in $\cF$ of any point in $\Gamma\backslash X$ is finite.

\subsection{Shimura varieties}\label{sec:shim-var}
We write $\bA_f$ for the ring of finite ad\`eles over $\bQ$.

Let $(\gG,\gX)$ be a Shimura datum and let $K$ be a compact open subgroup of $\gG(\bA_f)$. We obtain a double coset space
\[\Sh_K(\gG,\gX):=\gG(\bQ)\backslash[\gX\times(\gG(\bA_f)/K)],\]
where the implied action of $\gG(\bQ)$ is diagonal (the action on $\gX$ is given by conjugation and the action on $\gG(\bA_f)/K$ is given by left multiplication). Denote the equivalence class of $(x,g)\in\gX\times\gG(\bA_f)$ in $\Sh_K(\gG,\gX)$ by $[x,g]_K$ (or $[x,g]$ if there is no ambiguity).

The space $\Sh_K(\gG,\gX)$ is in natural bijection with a finite disjoint union of locally symmetric varieties. More precisely, if $X$ is a connected component of $\gX$ and $\cC$ is a (finite) set of representatives for $\gG(\bQ)_+\backslash\gG(\bA_f)/K$, then
\[\coprod_{g\in\cC}\Gamma_g\backslash X\to\Sh_K(\gG,\gX),\ \Gamma_gx\mapsto [x,g],\text{ where } \Gamma_g:=gKg^{-1}\cap\gG(\bQ)_+,\]
is a homeomorphism for the natural topologies
     \cite[Lem. 5.13]{Mil05}. Hence, $\Sh_K(\gG,\gX)$ is a (usually disconnected) quasi-projective complex algebraic variety, which we refer to as a {\it Shimura variety}. After possibly replacing $K$ by a finite index subgroup, it is smooth (any $K$ that is {\it neat}---see \cite[Sec. 4.1.4]{KY:AO}---would suffice).

\begin{remark}
    The language of Shimura varieties is clearly more involved than the language of locally symmetric varieties. Geometrically, the latter arguably gives a clearer picture, and will feature heavily when we discuss the Pila--Zannier strategy, where uniformisation maps $X\to\Gamma\backslash X$ play a key role due to the applicability of o-minimality. The language of Shimura varieties clarifies many arithmetic aspects. One example is the description of Hecke orbits (see Section \ref{sec:HCs}). Another is the following.
\end{remark}     

\subsubsection{Field of definition}

The complex variety $\Sh_K(\gG,\gX)$ has a {\it canonical model} over a number field $E(\gG,\gX)$ known as the {\it reflex field}. For the definition of canonical, see \cite[2.2.5]{Del79}. For the existence result, see \cite[Th. 7.2]{Milne1983}.

\begin{example}\label{ex:Y(1)}
    The simplest Shimura variety is the modular curve $Y(1):=\gSL_2(\bZ)\backslash\bH$, where $\bH$ denotes the upper half-plane. In the language of \ref{sec:shim-var}, this is $\Sh_{\gGL_2(\hat\bZ)}(\gGL_2,\gX)$, where $\hat\bZ:=\prod_p\bZ_p$ is the ring of integral ad\`eles and $\gX$ is the $\gGL_2(\bR)$-conjugacy class of the morphism $x_0:\bS\to\gGL_{2,\bR}$ given by 
    \[a+bi\mapsto\fullsmallmatrix{a}{b}{-b}{a}\]
    (which we identify with $\bH\sqcup\bar{\bH}$ via $x_0\mapsto i$).
    As an algebraic variety, $Y(1)\cong\bA^1$ (the affine line), and $E(\gGL_2,\gX)=\bQ$.
\end{example}

\begin{example}
    The most important Shimura variety (or, rather, {\it varieties}) is the moduli space $\cA_g$ of principally polarised abelian varieties of dimension $g$. This is the natural generalisation of Example \ref{ex:Y(1)}, in which we replace $\gGL_2$ with the general symplectic group $\gGSp_{2g}$ and $\bH$ with the Siegel upper half-space $\cH_g$. It has dimension $g(g+1)/2$ and the reflex field is also $\bQ$.

    The Baily-Borel compactification $\cA^{\rm BB}_g$ of $\cA_g=\gSp_{2g}(\bZ)\backslash\cH_g$ is stratified by locally closed subvarieties:
    \[\cA^{\rm BB}_g=\cA_g\sqcup\cA_{g-1}\sqcup\cdots\sqcup\cA_1\sqcup\cA_0,\]
    where by $\cA_0$ we denote a single point.
    \end{example}

\subsection{Generalisations}

Shimura varieties can be generalised in several directions. Relaxing the reductive condition on $\gG$ and replacing $\gX$ with the $\gG(\bR)\gU(\bC)$-conjugacy class of a homomorphism $\bS_\bC\to\gG_\bC$, for some normal unipotent subgroup $\gU$ of $\gG$, and imposing suitable conditions (see \cite[Def. 2.1]{pink:thesis}), we obtain via the same construction a more general class of algebraic varieties, known as {\it mixed Shimura varieties}, which include the universal families of abelian varieties.

In a different direction, by \cite[Sec. 2.3]{MoonenI}, for a suitable representation $\gG\to\gGL(V)$, with $V$ a finite-dimensional $\bQ$-vector space, and $K\subset\gG(\bA_f)$ sufficiently small, we obtain a so-called {\it variation of $\bQ$-Hodge structures} (in particular, a $\bQ$-local system with additional structure---see \cite[Sec. 1.2]{MoonenI}) on any connected component of $\Sh_K(\gG,\gX)$. This leads to the study of arbitrary smooth irreducible complex quasi-projective algebraic varieties endowed with such structures. This includes the local systems $R^i\!f_*\bQ$ for any smooth projective morphism $f:Y\to Z$ of smooth irreducible complex quasi-projective algebraic varieties. Via the associated {\it period map}, this setting is intimately related to the theory of arithmetic quotients $\Gamma\backslash\gG(\bR)/M$, where, now, $\gG$ is an arbitrary semisimple group over $\bQ$, $M$ is an arbitrary connected compact subgroup of $\gG(\bR)$, and $\Gamma$ is arithmetic. 

As we survey the results pertaining to Shimura varieties, we will endeavour to comment on the state of the art in the above levels of generality. We will see that our geometric understanding extends far into these settings, whereas many arithmetic questions remain wide open, even in the context of Shimura varieties.

Of course, the generalisations do not end here. Notably, the settings above are unified by the notion of so-called variations of {\it mixed Hodge structures}. There are also the {\it non-arithmetic} locally symmetric spaces, or the so-called {\it distinguished categories} of Barroero--Dill \cite{BD:dist}. However, in order not to stray too far from our primary subject matter, we will not comment much on the (substantial) progress reaching into these domains.

\subsection{Morphisms of Shimura data}
Let $(\gH,\gX_\gH)$ and $(\gG,\gX)$ be Shimura data. A morphism $(\gH,\gX_\gH)\to(\gG,\gX)$ of Shimura data is a homomorphism $\gH\to\gG$ which induces a map $\gX_\gH\to\gX$. If $(\gH,\gX_\gH)\to(\gG,\gX)$ is a morphism of Shimura data with $\gH$ a subgroup of $\gG$ and $\gH\to\gG$ the inclusion, we refer to $(\gH,\gX_\gH)$ as a Shimura subdatum of $(\gG,\gX)$. 

\subsection{Morphisms of Shimura varieties}
Let $(\gH,\gX_\gH)\to(\gG,\gX)$ be a morphism of Shimura data given by $\varphi:\gH\to\gG$. Let $K_\gH\subset\gH(\bA_f)$ and $K\subset\gG(\bA_f)$ be compact open subgroups such that $\varphi(K_\gH)\subset K$. By a theorem of Borel \cite{Borel72}, the natural map
\[f:\Sh_{K_\gH}(\gH,\gX_\gH)\to\Sh_K(\gG,\gX),\]
which is holomorphic on connected components,
extends to the Baily--Borel compactifications and, hence, is a morphism of algebraic varieties, which we refer to as a {\it morphism of Shimura varieties}. If $\ker(\varphi)^\circ$ is a torus, then $f$ is finite (cf. \cite[Fact 2.6(a)]{Pin05}). It is defined over the compositum $F:=E(\gH,\gX_\gH)\cdot E(\gG,\gX)$ of reflex fields \cite[Cor. 5.4]{Del71}, and we denote by $[f]$ the induced map on algebraic cycles. 

\subsection{Hecke correspondences}\label{sec:HCs}

Let $\Sh_K(\gG,\gX)$ be a Shimura variety and let $a\in\gG(\bA_f)$. There is a tautological isomorphism
\[\Sh_K(\gG,\gX)\to\Sh_{a^{-1}Ka}(\gG,\gX),\ [x,g]_K\mapsto [x,ga]_{a^{-1}Ka}.\]
We obtain a finite correspondence
\[\Sh_K(\gG,\gX)\leftarrow\Sh_{K\cap aKa^{-1}}(\gG,\gX)\to\Sh_{a^{-1}Ka\cap K}(\gG,\gX)\to\Sh_K(\gG,\gX),\]
where the middle arrow is the tautological isomorphism and the outer arrows are the finite morphisms induced by the identity map $\gG\to\gG$. We refer to such a correspondence as a {\it Hecke correspondence}, and we denote the map of algebraic cycles on $\Sh_K(\gG,\gX)$ by $T_{K,a}$ (or $T_a$ if there is no ambiguity).

\subsubsection{Hecke orbits}

Let $s\in\Sh_K(\gG,\gX)$. We refer to the set of points in the support of the $T_a(s)$ over all $a\in\gG(\bA_f)$
as the {\it Hecke orbit} of $s$. For $s\in\cA_g$, the Hecke orbit of $s$ corresponds to the (isomorphism classes of) abelian varieties possessing a {\it polarised} isogeny to the abelian variety associated with $s$ (see \cite[Sec. 1]{Orr15}).

Richard--Yafaev have generalised the notion of Hecke orbits in several ways \cite[Sec. 2]{RY:heights}. For $x\in\gX$, which we may consider as a morphism $\bS\to\gMT(x)_\bR$, they define the {\it generalised Hecke orbit} of $x$ in $X$ as the set
\[\cH(x):=\gX\cap\{\phi\circ x: \phi\in{\rm Hom}(\gMT(x),\gG)\}.\]
Then, for $[x,g]\in\Sh_K(\gG,\gX)$, they define the generalised Hecke orbit of $[x,g]$ as the image of $\cH(x)\times\gG(\bA_f)$ in $\Sh_K(\gG,\gX)$. This notion has several functoriality properties not enjoyed by classical Hecke orbits.

\subsection{Special subvarieties}
For any morphism of Shimura varieties $f:\Sh_{K_\gH}(\gH,\gX_\gH)\to\Sh_K(\gG,\gX)$, and any $g\in\gG(\bA_f)$, we refer to an irreducible component of the (Zariski closed) image of $T_g\circ [f]$ as a {\it special subvariety} of $\Sh_K(\gG,\gX)$. Note that the definition is unaffected if we assume that $\gH:=\gMT(\gX_\gH)$ or if we restrict to Shimura subdata $(\gH,\gX_\gH)$ of $(\gG,\gX)$. Special subvarieties of dimension zero are referred to as special points. 

\subsubsection{Examples} The special subvarieties of $Y(1)^n$ are defined by taking irreducible components of the loci defined by imposing the classical modular polynomial $\Phi_N(X,Y)=0$ on pairs of (not necessarily distinct) coordinates. A distinguished class of special subvarieties of $\cA_g$ (known as special subvarieties of {\it PEL type}) are defined by taking the irreducible components of the loci defined by $R\subset{\rm End}(A)$ for any ring $R$. In general, these are not the only special subvarieties of $\cA_g$, however (see \cite[Sec. 4]{Mum:PEL} for Mumford's famous example of non-PEL type special subvarieties of $\cA_4$).

\subsubsection{Defect} It is straightforward to show that the intersection of two special subvarieties of $S:=\Sh_K(\gG,\gX)$ is a finite union of special subvarieties. Therefore, for any irreducible subvariety $V$ of $S$ there is a smallest special subvariety of $S$ containing $V$, which we denote $\langle V\rangle$. Pink introduced the {\it defect}
\[\delta(V):=\dim\langle V\rangle-\dim V\]
of $V$ \cite{pink:generalisation}. If $\langle V\rangle$ is a connected component of $S$, we will follow convention and say that $V$ is {\it Hodge generic} in $S$.

\subsubsection{Complexity}

In order to count special subvarieties of $\Sh_K(\gG,\gX)$, and to compare arithmetic quantities associated with them, it is customary to attach to each special subvariety $Z$ a natural number $\Delta(Z)$ that captures (some of) its intrinsic {\it complexity}.

There are several natural ways to do this. For the special curve $Z$ in $Y(1)^2$ defined by $\Phi_N(X,Y)=0$, it is natural to define $\Delta(Z):=N$. On the other hand, if a special subvariety of $Y(1)^n$ has a fixed, necessarily special, coordinate (which arises when a condition of the form $\Phi_N(X,X)=0$ is imposed), then the absolute value of the discriminant of the associated CM elliptic curve is usually incorporated. More generally, for a PEL type special subvariety $Z$ of $\cA_g$ associated with a ring $R$, it is natural to consider $\Delta(Z)=|\disc(R)|$.

On the other hand, one can consider complexities defined in the $(\gG,\gX)$ language (as in \cite{UY:AO} and \cite{DO16}) or the underlying algebraic geometry (as in \cite[Def. 10.2]{DR18}). Urbanik has introduced a Hodge-theoretic complexity (see \cite{urbanik:degress}).

\subsubsection{Strongly and non-factor special subvarieties}

Let $Z$ be a special subvariety of $\Sh_K(\gG,\gX)$ associated with a Shimura subdatum $(\gH,\gX_\gH)$. We say that $Z$ is a {\it strongly special} subvariety if the image of $\gM:=\gMT(\gX_\gH)$ in $\gG^{\rm ad}$ is semisimple. 

More generally, following Ullmo \cite{ullmo:equidistribution}, we say that $Z$ is {\it non-factor} (or {\it non-facteur} or simply {\it NF}) if the image of $\gZ_\gG(\gM^\der)(\bR)$ in $\gG^{\ad}(\bR)$ is compact. (For the fact that a strongly special subvariety is non-factor, see \cite[Lem. 3.6]{DJK}.) This terminology is due to the fact that a special subvariety $Z$ of a connected Shimura variety $S$ is non-factor if and only if there exists no finite morphism of Shimura varieties
$S_1\times S_2 \to S$, with $\dim S_2>0$, such that $Z$ is equal to the image of some $S_1 \times \{z\}$ in $S$.

\subsection{Weakly special subvarieties}

Let $S:=\Sh_K(\gG,\gX)$ be a Shimura variety and let $(\gH,\gX_\gH)$ be a Shimura subdatum of $(\gG,\gX)$. Let $X_\gH$ be a connected component of $\gX_\gH$ and let $X_\gH=X_1\times X_2$ be a decomposition, associated with a decomposition $\gH^{\rm ad}=\gH_1\times\gH_2$. For any $x_2\in X_2$ and $g\in\gG(\bA_f)$, the image $Y$ of $X_1\times\{x_2\}\times\{g\}$ in $\Sh_K(\gG,\gX)$ is an irreducible (closed) subvariety, which we refer to as a {\it weakly special subvariety}. By \cite[Th. 4.3]{MoonenI}, an irreducible subvariety of $\Sh_K(\gG,\gX)$ is weakly special if and only if it is totally geodesic. Moreover, a weakly special subvariety of $S$ is a special subvariety of $S$ if and only if it contains a special point.

\subsubsection{Fibres of special subvarieties} Note that the image $Z$ of $X_\gH\times\{g\}$ in $S$ is a special subvariety of $S$ and we say that $Y$ is a {\it fibre} of $Z$.

\subsubsection{Weakly special defect} 

For any irreducible subvariety $V$ of $S$ there is a smallest weakly special subvariety of $S$ containing $V$, which we denote $\langle V\rangle_{\rm ws}$. We refer to
\[\delta_{\rm ws}(V):=\dim\langle V\rangle_{\rm ws}-\dim V\]
as the {\it weakly special defect} of $V$. Note that, for any irreducible subvarieties $W\subset V$ of $S$, 
\[\delta(V)-\delta_{\rm ws}(V)\leq\delta(W)-\delta_{\rm ws}(W)\]
\cite[Prop. 4.4]{DR18}. Habegger--Pila named this the {\it defect condition} \cite[Def. 4.2]{HP16}. Barroero--Dill show that the defect condition holds in any distinguished category \cite[Th. 7.2]{BD:dist}. In particular, this implies the condition for any connected mixed Shimura variety, which was shown independently by Cassani \cite{cassani:thesis}.

\subsection{Semi-algebraic sets}
A set $A\subset\bR^n$ is {\it semi-algebraic} if it is a finite union of sets defined by real polynomial equalities and inequalities. For an arbitrary set $A\subset\bR^n$ we denote by $A^{\rm alg}$ the union of the connected positive-dimensional semi-algebraic sets contained in $A$.

\subsection{o-minimal structures}

By a {\it structure}, we refer to a sequence $(S_n)_{n\in\mathbb{N}}$, where, for each $n\in\bN$, $S_n$ is a Boolean algebra (under the standard set-theoretic operations) of subsets of $\bR^n$ which contains all of the semi-algebraic subsets, such that, for $n\geq m\in\bN$,
\begin{enumerate}
    \item if $A\in S_n$ and $B\in S_m$, then $A\times B\in S_{n+m}$;
    \item if $A\in S_n$ and $\pi:\bR^n\to\bR^m$ is the projection onto the first $m$ coordinates, then $\pi(A)\in S_m$.
    \end{enumerate}
It is a famous result of Tarski--Seidenberg that the semi-algebraic sets themselves form a structure \cite[Th. 8.6.6]{Mishra:TS}.
    
Given a structure $\cS=(S_n)_{n\in\bN}$, we say that a set $A\subset\bR^n$ is definable in $\cS$ if $A\in S_n$. We say that $\cS$ is {\it o-minimal} if, for every $A\in S_1$, the boundary $\partial A$ is finite. The o-minimal structure of greatest significance to us will be the smallest structure containing all sets defined by the real exponential function and the restricted analytic functions, denoted $\bR_{\rm an,exp}$. For the fact that this structure is indeed o-minimal (and related references), see \cite{vdDM:Ranexp}, wherein the authors also establish that $\bR_{\rm an,exp}$ has so-called {\it analytic cell decomposition}.

\subsection{Heights}
For a polynomial $f=a_kX^k+\cdots a_1X+a_0\in\bZ[X]$, we define $H(f)$ to be the maximum of $|a_0|,\ldots,|a_k|$. For $k\in\bN$ and $x\in\bR$, we define $H_k(x)$ to be the minimum of $H(f)$ as $f$ varies over all primitive integral polynomials such that $f(x)=0$ and $\deg(f)\leq k$. (If there is no such $f$, which is to say $[\bQ(x):\bQ]>k$, we define $H_k(x)=+\infty$.) For $x=(x_1,\ldots,x_n)\in\bR^n$, we define $H_k(x)$ to be the maximum of $H_k(x_1),\ldots,H_k(x_n)$.

For $x\in\bP^n(\Qbar)$, we define the multiplicative (resp. logarithmic) Weil height $H(x)$ (resp. $h(x)$) of $x$ as in \cite[Def. 1.5.4]{bombieri:Gfns}.

\subsection{The Pila--Wilkie theorem}

The following theorem due to Pila--Wilkie has revolutionised the field of unlikely intersections.

\begin{theorem}[Pila--Wilkie, {\cite{PW:counting,pila:counting}}]\label{teo:PW}
    Let $\cS$ be an o-minimal structure and let $A\subset\bR^n$ be definable in $\cS$. Let $k\in\bN$ and let $\epsilon>0$. Then there exists $c:=c(A,k,\epsilon)>0$ such that, for any $T\geq 1$,
    \[|\{x\in A\setminus A^{\rm alg}:H_k(x)\leq T\}|\leq cT^\epsilon.\]
\end{theorem}

The Pila--Wilkie theorem has a long history, but we will only mention the seminal 1989 paper of Bombieri--Pila \cite{BP}. We also emphasise, however, that generalising the theorem is an important area of contemporary research. Indeed, for our purposes, we will need a modest strengthening capable of handling so-called {\it semi-rational} points. For this, we refer to a simple case of the versions appearing in \cite[Sec. 7]{HP16} (cf. \cite[Th. 9.1]{DR18}).

\begin{theorem}[Semi-rational Pila--Wilkie]\label{teo:PWSR}
    Let $\cS$ be an o-minimal structure admitting analytic cell decomposition and let $\cD\subset\bR^m\times\bR^n=\bR^{m+n}$ be definable in $\cS$. Let $k\in\bN$ and let $\epsilon>0$. Then there exists $c:=c(\cD,k,\epsilon)>0$ such that, for any $T\geq 1$, if
    \[\Lambda\subset\{y\in\bR^n:\exists x\in\Qbar^m\text{ with }H_k(x)\leq T\text{ and }(x,y)\in\cD\}\]
    satisfies $|\Lambda|>cT^\epsilon$, there exists a continuous definable function $\beta:[0,1]\to\cD$, which is real analytic on $(0,1)$, such that $\pi_1\circ\beta$ (resp. $\pi_2\circ\beta$) is semi-algebraic (resp. non-constant), and satisfies $\pi_2(\beta(0))\in\Lambda$, where \[\pi_1:\bR^m\times\bR^n\to\bR^m\text{ and }\pi_2:\bR^m\times\bR^n\to\bR^n\] 
    denote the projections to the first and second factors.
\end{theorem}

\section{A review of problems of unlikely intersections}

\subsection{Zilber--Pink} The central problem of unlikely intersections in Shimura varieties is the {\it Zilber--Pink conjecture}. Pink first formulated this conjecture for (mixed) Shimura varieties in \cite[Conj. 1.3]{pink:generalisation}, whereas Zilber had previously formulated an analogous conjecture for semiabelian varieties in \cite[Conj. 2]{Zilber02}. Bombieri--Masser--Zannier first formulated and proved a case of Zilber's conjecture for curves in algebraic tori \cite[Th. 2]{BMZ99}, at the same time raising the analogous question for abelian varieties. 

There was a diversity in perspectives that led these authors to their conjectures, but questions of transcendence played a key role:
Bombieri--Masser--Zannier were interested in the multiplicative independence of algebraic numbers, and Zilber was interested in Schanuel's conjecture. Andr\'e, for whom the Andr\'e--Oort and Andr\'e--Pink--Zannier conjectures (see below) are named, was, in part, interested in Grothendieck's period conjecture. 

\begin{conjecture}[Zilber--Pink]\label{ZP}
    Let $S$ be a Shimura variety and let $V$ be an irreducible subvariety of $S$. Suppose that the intersection
of $V$ with the union of all special subvarieties of $S$ of codimension greater than $\dim(V)$ is Zariski dense in $V$. Then $V$ is contained in a proper special subvariety of $S$.
\end{conjecture}

Stated this way, the import of the term {\it unlikely intersections} is self-evident; randomly chosen subvarieties $V$ and $Z$ satisfying ${\rm codim}(Z)>\dim(V)$ are {\it unlikely} to produce a non-empty intersection.

\subsubsection{Other formulations}

Zilber--Pink has several other formulations. The following definition imitates the terminology Zilber used to make his conjecture for semiabelian varieties.

\begin{definition}
    Let $S$ be a Shimura variety and let $V$ be an irreducible subvariety. An irreducible subvariety $W$ of $V$ is an {\it atypical subvariety} of $V$ if there is a special subvariety $Z$ of $S$ such that
    \begin{enumerate}
        \item $W$ is an irreducible component of $V\cap Z$;
        \item $\dim W>\dim V+\dim Z-\dim S$.
    \end{enumerate}
    We denote by $\Atyp(V)$ the union of the atypical subvarieties of $V$.
\end{definition}

An ostensibly stronger formulation of Zilber--Pink can be stated as follows. (For a proof of Conjecture \ref{ZP} for a given $S$ assuming Conjecture \ref{ZP:atyp} for $S$ see \cite[Lem. 3.6]{DR18}.)

\begin{conjecture}[Strong Zilber--Pink I, {\cite[Conj. 2.2]{HP16}}]\label{ZP:atyp}
    Let $S$ be a Shimura variety and let $V$ be an irreducible subvariety of $S$. Then $\Atyp(V)$ is equal to a finite union of atypical subvarieties of $V$.
\end{conjecture}

Habegger--Pila gave another formulation using the following terminology (see \cite[Def. 2.5]{HP16}).

\begin{definition}
    Let $S$ be a Shimura variety and let $V$ be an irreducible subvariety. An irreducible subvariety $W$ of $V$ is an {\it optimal subvariety} of $V$ if for any irreducible subvariety $Y$ of $V$ strictly containing $W$, we have $\delta(Y)>\delta(W)$. An {\it optimal point} is an optimal subvariety of dimension zero.
    
    We denote by $\Opt(V)$ the set of optimal subvarieties of $V$ and by $\Opt_0(V)$ the set of optimal points of $V$.
\end{definition}

The formulation of Habegger--Pila is as follows. (For a proof of Conjecture \ref{ZP:atyp} for a given $S$ assuming Conjecture \ref{ZP:opt} for $S$, and vice versa, see \cite[Lem. 2.7]{HP16}.)

\begin{conjecture}[Strong Zilber--Pink II, {\cite[Conj. 2.6]{HP16}}]\label{ZP:opt}
    Let $S$ be a Shimura variety and let $V$ be an irreducible subvariety of $S$. Then $\Opt(V)$ is finite.
\end{conjecture}

The name {\it Strong Zilber--Pink} is in some sense a misnomer, however. Barroero--Dill have shown that Conjecture \ref{ZP:opt} is true (that is, for all $S$) if and only if Conjecture \ref{ZP} is true (for all $S$) \cite[Th. 12.4]{BD:dist}.

\subsubsection{Generalisations}

More recently, Klingler formulated a Zilber--Pink conjecture for any variation of mixed Hodge structures, vastly generalising the conjectures above \cite[Sec. 1.5]{Klingler:ZP}. Baldi--Klingler--Ullmo further strengthened Klingler's conjectures for variations of Hodge structures \cite[Conj. 2.5]{BKU}. At the time of writing, the most general formulation, for variations of mixed Hodge structures, is \cite[Conj. 4.9]{BU:effective}. 

\subsection{Andr\'e--Oort}

The Zilber--Pink conjecture has many notable corollaries. The earliest and most prominent of these is the {\it Andr\'e--Oort conjecture}. This statement concerns the distribution of special points and combines conjectures of Andr\'e \cite[Ch. X, Sec. 4.3, Prob. 1]{And89} and Oort \cite{Oort94}. Its proof will be discussed in Section \ref{subsec:AO}.

\begin{theorem}[Andr\'e--Oort, {\cite[Th. 1.1]{PST+}}]
    Let $S$ be a Shimura variety and let $V$ be an irreducible subvariety of $S$. Suppose that the special points of $S$ contained in $V$ are Zariski dense in $V$. Then $V$ is a special subvariety of $S$.
\end{theorem}

Note that the Andr\'e--Oort conjecture is the formal analogue of the Manin--Mumford conjecture for abelian varieties and its counterpart for algebraic tori, which originated in conjectures of Lang (see \cite{Zan12} for more on the history of these statements).

\subsection{Andr\'e--Pink--Zannier}

Another notable corollary of Zilber--Pink is the so-called {\it Andr\'e--Pink--Zannier conjecture}, concerning the distribution of points in a single Hecke orbit. It was first formulated for curves in \cite[Ch. X, Sec. 4.5, Prob. 3]{And89} and for arbitrary mixed Shimura varieties in \cite[Conj. 1.6]{Pin05}. It is mentioned in passing in \cite[p. 12]{yaf:thesis}, and it is closely related to a conjecture attributed to Zannier (see \cite[Conj. 1.4]{Geo:APZ}).

\begin{conjecture}[Andr\'e--Pink--Zannier]\label{teo:APZ}
    Let $S$ be a Shimura variety and let $V$ be an irreducible subvariety. Suppose that the intersection of $V$ with a single Hecke orbit in $S$ is Zariski dense in $V$. Then $V$ is a weakly special subvariety of $S$.
\end{conjecture}

Notably, Conjecture \ref{teo:APZ} has recently been proved for Shimura varieties of abelian type by Richard--Yafaev \cite{RY:APZ} (see \cite[p. 4118]{ullmo:points} for the precise definition of Shimura varieties of abelian type). In fact, they prove the result for generalised Hecke orbits. 

\subsection{Mordell--Lang for Shimura varieties}

In \cite{HP12}, Habegger--Pila established an analogue of the Mordell--Lang conjecture for products of modular curves. In \cite{AD}, we observed that the generalisation to arbitrary Shimura varieties is equivalent to the conjunction of Andr\'e--Oort and 
Andr\'e--Pink--Zannier. Richard--Yafaev have formulated a stronger conjecture \cite[Conj. 5.2]{RY:generalisation} and have obtained this conjecture for Shimura varieties of abelian type \cite[Th. 1.4]{RY:hybrid}. As a consequence, they deduce Zilber--Pink for Hodge generic weakly special subvarieties of codimension $1$.

For a proof that Andr\'e--Pink--Zannier (and, indeed, the aforementioned strengthening of Richard--Yafaev) is a corollary of Zilber--Pink, see \cite[Sec. 7]{RY:generalisation}. 

\subsection{Finite characteristic}

An Andr\'e--Oort conjecture for $\bA^2_{\bF_p}$ was formulated and proved under the generalised Riemann hypothesis (GRH) for quadratic fields by Edixhoven-Richard \cite{ER:modp}. This inspired an Andr\'e--Oort conjecture over the ring of integers of a number field \cite{Richard:AO,Saetonne} (see
\cite{BRU:MM} for an analogue in the abelian setting). Shankar--Lam formulate an Andr\'e--Pink--Zannier conjecture in characteristic $p$ and prove a special case \cite{LS:APZ}.

\section{A history of the results}

\subsection{Andr\'e--Oort}\label{sec:AOhistory}

In 1998, Moonen proved Andr\'e--Oort in $\cA_g$ under a mod $p$ assumption \cite{MoonenII}, and Yafaev generalised this to arbitrary Shimura varieties shortly thereafter \cite{Yafaev:Moonen}. However, the first strategy towards a full proof of Andr\'e--Oort originated in a paper of Edixhoven \cite{Edixhoven:2prod}, wherein the author proved the conjecture for $S=Y(1)^2$ assuming GRH for imaginary quadratic fields. Edixhoven subsequently applied his approach, still under GRH, to arbitrary products of modular curves \cite{Edixhoven:nprod} and Hilbert modular surfaces \cite{Edixhoven:hilbert}. Edixhoven--Yafaev extended the strategy to curves in arbitrary Shimura varieties containing infinitely many special points in a Hecke orbit, and Yafaev generalised this result to arbitrary collections of special points under GRH \cite{yafaev:AO}. Another decade of work culminated in a proof, due to Klingler--Ullmo--Yafaev, of Andr\'e--Oort for arbitrary Shimura varieties assuming the generalised Riemann hypothesis for CM fields (GRHCM) \cite{UY:AO,KY:AO}. This proof will be summarised in Section \ref{sec:EKUY}.

Except for Andr\'e's proof of Andr\'e--Oort in $Y(1)^2$ \cite{And98}, unconditional results were lacking until, in 2008, a new strategy towards problems of unlikely intersections emerged, first in a new proof, due to Pila--Zannier, of the Manin--Mumford conjecture \cite{PZ:manin-mumford}, and realised shortly thereafter in the setting of Shimura varieties by Pila in a proof of Andr\'e--Oort for $Y(1)^n$ \cite{pila:y(1)n,pila:y(1)2}. This inspired unconditional proofs for Hilbert modular surfaces \cite{dy:hilbert}, $\cA_2$ \cite{PT13}, 
$\cA_g$ \cite{Tsim:AOAg} and, eventually, general Shimura varieties \cite{PST+}, though these results required several other significant advances. The Andr\'e--Oort conjecture for mixed Shimura varieties follows from a theorem of Gao \cite{Gao:reduction}. This strategy, known as the Pila--Zannier strategy, will be covered in Section \ref{sec:PZ}.

Effective results (that is, quantitative results whose outputs are shown to be computable) remain relatively sparse. K\"uhne \cite{Kuh12} and Bilu--Masser--Zannier independently gave an effective proof of Andr\'e--Oort for $Y(1)^2$. Binyamini has given some effective results for $Y(1)^n$ \cite{Binyamini:AOY1}, whereas Breuer gave an effective proof for curves in $Y(1)^n$ under GRH for imaginary quadratic fields \cite{Breuer:AO}. 

More general effective results may be on the horizon. Binyamini has recently announced an effective variant of the Pila--Wilkie theorem that is applicable in all of the settings hitherto described \cite{Binyamini:Log}. Another aspect of modern techniques, based on an idea of Andr\'e using {\it G-functions}, which will feature heavily below, is also effective, and was employed by Binyamini--Masser to yield effective results in Hilbert modular varieties \cite{BM:AO}.

\subsubsection{Applications}

Oort was led to formulate what we now refer to as the Andr\'e--Oort conjecture for $\cA_g$ in the knowledge that it reduced a question of Coleman on CM Jacobians to the question of whether the so-called Torelli locus $\cT_g\subset\cA_g$ contains any positive-dimensional special subvarieties meeting the image of the Torelli morphism $\cM_g\to\cA_g$. For the current state of Coleman's fascinating question, see \cite{Moonen:Coleman}.

Edixhoven--Yafaev were prompted to prove their aforementioned result following a suggestion of Cohen--W\"ustholz \cite{CW:hyper}, in order to conclude a programme of Wolfart on the algebraicity of values of hypergeometric functions at algebraic numbers \cite{Wolfart}. 

Chai--Oort showed that the Andr\'e--Oort conjecture gives an affirmative answer to the question of Katz--Oort as to whether there exists an abelian variety over $\Qbar$ which is not isogenous to a Jacobian \cite{CO:Jac}. Tsimerman modified the argument to give an unconditional answer in \cite{Tsim:Jac}. Recently, Masser--Zannier gave several strong refinements using the Pila--Zannier method \cite{MZ:Jac}. Shankar--Tsimerman have investigated the question in characteristic $p$ \cite{ST:Jac1,ST:Jac2}.

In \cite[Th. 1.3]{AD}, we use Andr\'e--Oort to derive a special case of the Zilber--Pink conjecture.

\subsection{Andr\'e--Pink--Zannier}

The aforementioned result of Edixhoven--Yafaev \cite{EY:AO} can be interpreted as Andr\'e--Pink--Zannier for curves intersecting Hecke orbits of special points.
In \cite[Th. 7.6]{Pin05}, Pink proved Andr\'e--Pink--Zannier for Hecke orbits of so-called {\it Galois generic points} in $\cA_g$, appealing to the equidistribution of Hecke points proved by Clozel--Oh--Ullmo \cite{COU:equidistribution}, and Cadoret--Kret generalised this result to arbitrary Shimura varieties defined by $\gG$ almost $\bQ$--simple in \cite[Th. B]{CK:APZ}, also using equidistribution. 

Once more, the Pila--Zannier strategy has paved the way towards more general results. Orr applied it to prove Andr\'e--Pink--Zannier for curves in $\cA_g$ \cite{Orr15}, and it forms the basis for the aforementioned work \cite{RY:heights,RY:APZ} of Richard--Yafaev. It has also been used to prove several special cases of Andr\'e--Pink--Zannier in mixed Shimura varieties (see \cite[Sec. 1]{Dil}).

\subsection{Zilber--Pink}\label{sec:ZPhistory}

For varieties not definable over $\Qbar$ much is now known. See \cite[Cor. 1.7]{KT:notQbar} for the most general result in arbitrary Shimura varieties, and \cite[Th. 1.1]{BD:Hecke} for Hodge generic subvarieties of $\cA_g$. Pila previously handled curves in $Y(1)^3$ \cite[Th. 1.4]{pila:fermat}, and Aslanyan--Eterovic--Fowler prove a generalisation (not implied by the aforementioned results) in $Y(1)^n$ \cite[Th. 1.1]{AEF:Qbar}.

Beyond Andr\'e--Oort and Andr\'e--Pink--Zannier, however, cases of Zilber--Pink over $\Qbar$ are relatively sparse. Habegger--Pila proved Zilber--Pink for curves in $Y(1)^n$ that are {\it asymmetric} \cite{HP12} and Orr generalised their arguments to asymmetric curves in $\cA^2_g$. 

In \cite{ExCM,QRTUI,LGO}, we obtained Zilber--Pink for curves in $\cA_2$ whose Zariski closures in $\cA^{\rm BB}_2$ intersect $\cA_0$. Our works \cite{PEL,LGO} combined with the work of Bhatta \cite{Bhatta}, yield ``simple PEL type'' Zilber--Pink for Hodge generic curves in $\cA_g$ whose Zariski closures in $\cA^{\rm BB}_g$ intersect $\cA_0$. Papas establishes that such curves, when contained in the Torelli locus, contain only finitely many points corresponding to non-simple Jacobians \cite[Th. 1.1]{PapasTg}. In \cite[Th. 1.4]{PapasAg}, Papas obtains cases of simple PEL type Zilber--Pink for Hodge generic curves in $\cA_g$ whose Zariski closures in $\cA^{\rm BB}_g$ intersect $\cA_{g-h}$ with $h\geq 2$.

In \cite{Y(1)}, we proved Zilber--Pink for curves in $Y(1)^n$ whose Zariski closures contain $(\infty,\ldots,\infty)$. Papas extends this to any curve in $Y(1)^3$ whose Zariski closure contains a special point in the boundary \cite[Th. 1.5]{PapasY(1)}. In \cite{DOP}, we handle the case of a curve in $Y(1)^3$ whose Zariski closure contains a point in the boundary lying on a modular curve. 

\subsubsection{Other settings}

Zilber--Pink for curves in $\bG^n_m$ defined over $\Qbar$ was proved by Maurin \cite{Maurin:ZP}. The analogous statement for abelian varieties was obtained by Habegger--Pila \cite{HP16} and, for semiabelian varieties, it is due to Barroero--K\"uhne--Schmidt \cite{BKS:ZP}. Both of these proofs use the Pila--Zannier method. Maurin's theorem was extended to curves over $\bC$ by Bombieri--Masser--Zannier \cite{BMZ08}, and the extension for abelian varieties is obtained by Barroero--Dill in \cite{BD:abelian}. For semiabelian varieties, see \cite[Sec. 14]{BD:dist}. 

Zilber--Pink in the mixed setting has a rich history. Earlier works, due to Masser--Zannier and Corvaja--Masser--Zannier, focused on the so-called relative Manin-Mumford conjecture \cite{MZ:2008,MZ:2010,MZ:2012,MZ:2014,MZ:2015,CMZ:RMM}. These culminated in a proof of relative Manin--Mumford for curves in abelian schemes defined over $\Qbar$ and facilitated striking applications to the solvability of Pell's equation and to integration in elementary terms \cite{MZ:integration}. Recently, Gao--Habegger proved the relative Manin--Mumford conjecture for any subvariety of an abelian scheme (in characteristic $0$) \cite{GH:RMM}. 

A striking observation of Bertrand was that relative Manin--Mumford is false for semiabelian schemes \cite{Bertrand:counter}! The correct statement (implied by Zilber--Pink) for one-dimensional families of semi-abelian surfaces was obtained by Bertrand--Masser--Pillay--Zannier \cite{BMPZ:RRM}. For intersections between curves and subgroup schemes, see the work of Barroero--Capuano \cite{BC:ab-schemes} and the recent survey of Capuano \cite{Capuano:survey}. 

Note that all of the aforementioned works in the relative setting invoke the Pila--Zannier strategy. Another important development in the relative setting concerns the so-called {\it Betti map}, which locally identifies the fibres of an abelian scheme of relative dimension $g$ with the real torus $(\bR/\bZ)^{2g}$. Gao's results on its generic rank along a subvariety and the associated {\it degeneracy loci} \cite{Gao:Betti} have played an important role in recent developments in Diophantine geometry, including the uniform Mordell--Lang conjecture \cite{DGH:UMLC,Kuh:UML,GGK:UML}, Ge's extension to abelian schemes of Habegger's bounded height theorem \cite{Ge:bounded-height}, and the relative Manin--Mumford conjecture mentioned above.

\section{Andr\'e--Oort: the Edixhoven--Klingler--Ullmo--Yafaev strategy }\label{sec:EKUY}

In this section, we loosely summarise the first general strategy towards Andr\'e--Oort. As mentioned above, its first incarnation appeared in a paper of Edixhoven, but it clearly owes much to the approach of Hindry in the abelian setting \cite{Hindry}. It is interesting to note that some of these ideas have recently been revived in finite characteristic. Be aware that our notations $V$ and $Z$ are inverted in \cite{KY:AO}.
 
\subsection{Overview} Suppose that $S$ is defined by a Shimura datum $(\gG,\gX)$. Using monodromy arguments, Klingler--Yafaev show that, given 
\begin{enumerate}
    \item a special but non-strongly special subvariety $Z$ of $S$ contained in $V$;
    \item a suitable $m\in\gG(\bQ_p)$, for a suitable prime $p$,
\end{enumerate}
the inclusion $V\subset T_m(V)$ implies that $V$ contains a special subvariety that contains $Z$ properly. The idea is to then apply this fact iteratively to a (Zariski dense) set $\cS$ of special subvarieties (initially of dimension $0$) in $V$ until $V\in\cS$ is forced, simply by dimension.

\subsection{Obtaining the ingredients} The existence of suitable $m\in\gG(\bQ_p)$ is established using {\it Bruhat--Tits buildings}. This requires that $p$ splits in the splitting field of the torus $\gZ:=\gZ(\gH)^{\circ}$ with $(\gH,\gX_\gH)=(\gMT(\gX_\gH),\gX_\gH)$ the Shimura subdatum defining $Z$. 

To obtain $V\subset T_m(V)$, the idea is as follows. Suppose that $r:=\dim V-\dim Z=1$ (the general case proceeds by induction on $r$). Ullmo--Yafaev show under GRHCM that the degree of the Galois orbit $Z^{\rm Gal}$ of $Z$ is bounded from below in terms of a complexity $\Delta_{\rm UY}(Z)$ (which depends only on $\gZ$). Moreover, after some technical steps, $Z^{\rm Gal}$ is contained in $V\cap T_m(V)$ and so, if $\Delta_{\rm UY}(Z)$ suitably exceeds the degree of $V\cap T_m(V)$, the latter cannot be a proper intersection and the desired conclusion is forced. Up to a constant, the degree of $V\cap T_m(V)$ is the degree of $T_m$, and Klingler--Yafaev show that there exist suitable $m$ for which this is bounded from above by a uniform power of $p$. In this way, bearing in mind the requirement on $p$ above, the problem is reduced to the availability of small split primes, hence GRHCM (again).

\subsection{Overcoming the obstacles} 
Clearly, the strategy breaks down if the iterative process produces only strongly special subvarieties, or, more generally, if $\Delta_{\rm UY}(Z)$ is bounded for $Z\in\cS$. Ullmo--Yafaev show that, essentially, the latter can only happen for sets of strongly special subvarieties. Therefore, it remains to show that the Zariski closure of a union of strongly special subvarieties is a finite union of special subvarieties.

There are at least two viable approaches to this problem. The proof of Klingler--Ullmo--Yafaev appeals to the equidistribution of strongly special subvarieties proved by Clozel--Ullmo \cite{CU:equidistribution}. Using Prasad's volume formula \cite{Prasad:volumes}, the present author obtained lower bounds for the degrees of strongly special subvarieties and used these to extend the above strategy to all special subvarieties \cite{daw:degrees}. Both of these proofs are unconditional and have been generalised to non-factor special subvarieties \cite{ullmo:equidistribution,DJK}. The latter is also effective.

\begin{remark}\label{KUY:unconditional}
    In Valentin Blomer's recent Number Theory Web Seminar \cite{Blomer}, Farrell Brumley raised the question of whether the Quasi-Riemann hypothesis, announced by OpenAI in \cite{OAI}, and its extension to Dedekind zeta functions by Cureton--Thorner \cite{hexagon:2610.00200v1}, can replace the use of GRHCM in this proof. This is indeed the case, due to recent work of Frei--Sofos--Widmer, mentioned by Brumley, which proves a sufficiently strong effective Chebotarev estimate \cite[Corollary 1.2]{FSW}, which is unconditional given the theorem of Cureton--Thorner. The use of GRHCM above is, in both instances, to facilitate an effective Chebotarev estimate of precisely this form in CM fields. We are most grateful to Philipp Habegger for bringing this to our attention. An effective implementation of this has been announced in \cite{jaswal_2026_24337359}.
\end{remark}

\section{The Pila--Zannier strategy}\label{sec:PZ}

The starting point for the Pila--Zannier strategy is the pivotal observation that, for any locally symmetric variety $S=\Gamma\backslash X$, the uniformisation map $\pi:X\to S$, while transcendental for the algebraic structures on $X$ and $S$, is definable in $\bR_{\rm an,exp}$ when restricted to a fundamental set $\cF$ for $\Gamma$. This fact was first obtained for $S=\cA_g$ by Peterzil--Starchenko \cite{PS:definability} and, then, for general Shimura varieties by Klingler--Ullmo--Yafaev \cite{kuy:ax-lindemann}. Moreover, Gao observed that the latter yields the result for mixed Shimura varieties \cite[Sec. 10.1]{Gao:AxL}, and it was generalised to arithmetic quotients by Bakker--Klingler--Tsimerman \cite[Th. 1.1(1)]{BKT:definability}. For variations of mixed Hodge structures, see \cite{BBKT:definability}. 

In particular, it follows that, for any algebraic subvariety $V$ of $S$, the set $\cV:=\pi^{-1}(V)\cap\cF$ is definable in $\bR_{\rm an,exp}$.

\subsection{Andr\'e--Oort}\label{subsec:AO}

It is useful to describe the Pila--Zannier strategy as having two components: one geometric, one arithmetic. The geometric component for Andr\'e--Oort is the following result due to Ullmo (reformulated slightly).

\begin{theorem}[Geometric Andr\'e--Oort, {\cite{Ull14}}]\label{teo:GAO}
Let $S$ be a Shimura variety and let $V$ be an irreducible subvariety of $S$. Then there exists a finite set $\Theta$ of special subvarieties of $S$ such that each maximal weakly special subvariety $W$ of $V$ is a fibre of some $Z\in\Theta$.
\end{theorem}

Alongside o-minimality, the main ingredient in the proof of Theorem \ref{teo:GAO} is the so-called {\it Ax--Lindemann} (or {\it Ax--Lindemann--Weierstrass}) {\it theorem} for Shimura varieties, an analogue of the result (without derivatives) of Ax for the exponential function \cite{Ax}, itself a functional analogue of the classical transcendence result due to Lindemann--Weierstrass. In the generality below, it was proved by Klingler--Ullmo--Yafaev \cite{kuy:ax-lindemann}, following earlier proofs for $\cA_g$ \cite{PT:AxL} and compact Shimura varieties \cite{UY:AxL}. It was proved for mixed Shimura varieties by Gao \cite{Gao:AxL}. All of these proofs rely on o-minimality.

\begin{theorem}[Ax--Lindemann for Shimura varieties]\label{teo:ALW}
Let $S=\Gamma\backslash X$ be a connected component of a Shimura variety and let $\pi:X\to S$. Let $V$ be an irreducible subvariety of $S$ and let $W$ be an irreducible algebraic subvariety of $X$ (see Section \ref{sec:alg-structure}) contained in $\pi^{-1}(V)$ and maximal for this property. Then $\pi(W)$ is a weakly special subvariety of $S$.
\end{theorem}

More recently, Richard--Ullmo reproved Ullmo's result \cite{RU:AO} using the results of Gorodnik, Ullmo and the present author on sequences of homogeneous probability measures on arithmetic locally symmetric spaces \cite{DGU,DGU2}. In fact, they prove a dynamical result in the generality of arithmetic quotients and apply this, in collaboration with Chen, to obtain geometric Andr\'e--Oort for variations of Hodge structures.  

\quad

It follows easily from geometric Andr\'e--Oort that, to prove the Andr\'e--Oort conjecture, it suffices to prove the following arithmetic statement.

\begin{theorem}[Arithmetic Andr\'e-Oort]\label{teo:AAO}
Let $S$ be a Shimura variety and let $V$ be an irreducible subvariety of $S$. Then $V$ contains only finitely many special points not contained in a positive-dimensional special subvariety contained in $V$.
\end{theorem}

The proof of Theorem \ref{teo:AAO} requires two arithmetic ingredients. In order to state them, let $\Sigma_0$ denote the set of special points of $S$ and fix a complexity $\Delta:\Sigma_0\to\bN$
such that, for any $B>0$, the set
\[\{s\in\Sigma_0:\Delta(s)<B\}\]
is finite. Replace $S$ with the connected component $\Gamma\backslash X$ containing $V$ and let $E$ be a number field over which $S$ and $V$ are defined. (Indeed, since special points are defined over $\Qbar$, we may assume that so too is $V$.)  Let $\pi:X\to S$ be the uniformisation map and fix a fundamental set $\cF$ in $X$ for $\Gamma$.

\subsubsection{Large Galois orbits for special points}\label{sec:galois-sp-pts}
The primary arithmetic ingredient is to show that there exists $\newC{LGO0-exp}>0$ such that, for any $s\in\Sigma_0$, we have 
    \[|\Gal(\Qbar/E)\cdot s|\gg\Delta(s)^{\refC{LGO0-exp}}.\]
This was originally conjectured for $S=\cA_g$ by Edixhoven \cite[Prob. 14]{EM:problems}, inspired by his result for Hilbert modular surfaces \cite{Edixhoven:hilbert}. The case $S=Y(1)$ follows immediately from the Brauer--Siegel theorem \cite{brauer:siegel}. Early techniques for general Shimura varieties were obtained in \cite{CU:tori}. The cases $S=\cA_g$ for $g\leq 6$ were obtained by Tsimerman in \cite{tsim:galois} (see also \cite{uy:galois} for $g\leq 3$), and Edixhoven's conjecture was later fully resolved by Tsimerman in \cite{Tsim:AOAg} using the Colmez conjecture on average, which had recently been obtained by Yuan--Zhang \cite{YZ:avColmez} and Andreatta--Goren--Howard--Madapusi-Pera \cite{AGHMP:avColmez}, independently of each other.

The proof for $S=\cA_g$ relied on Masser--W\"ustholz isogeny estimates \cite{MW95}. A new approach, applicable to general Shimura varieties, was obtained by Binyamini--Schmidt--Yafaev in \cite{BSY}. Large Galois orbits for special points have now been announced in full generality by Pila--Shankar--Tsimerman using new results of Esnault--Groechenig \cite{PST+}.

\subsubsection{Heights of pre-special points}\label{sec:hts-pre-sp}
We refer to a point $x\in X$ as pre-special if $\pi(x)\in\Sigma_0$. By \cite[Prop. 3.7]{UY:characterisation}, pre-special points belong to $X^\vee(\Qbar)$. In fact, as explained in \cite[Sec. 1.2]{DO16}, given an embedding $X^{\vee}\subset\bP^n_{\Qbar}$, there exists $k\in\bN$ such that the coordinates of any pre-special point are algebraic of degree at most $k$.
The second ingredient of the strategy is to show that there exists $\newC{Param0-exp}>0$ such that, for any pre-special $x\in\cF$, we have
\[{H}(x)\ll\Delta(\pi(x))^{\refC{Param0-exp}}.\]
For $S=Y(1)$, this is a relatively straightforward computation \cite[Sec. 4]{pila:counting}. For Hilbert modular surfaces, it appears in \cite[Th. 2.2]{dy:hilbert}. For $S=\cA_g$, it is \cite[Th. 3.1]{PT13}. For general Shimura varieties, it is the main result of \cite{DO16}.

\subsubsection{Proof of Theorem \ref{teo:AAO}}

Suppose $s\in V$ is a special point. Then, for each point $\sigma(s)\in\Gal(\Qbar/E)\cdot s$ and $x_\sigma\in\cV$ such that $\pi(x_\sigma)=\sigma(s)$, we have
\[{H}(x_\sigma)\ll\Delta(\sigma(s))^{\refC{Param0-exp}}\ll|\Gal(\Qbar/E)\cdot \sigma(s)|^{\frac{\refC{Param0-exp}}{\refC{LGO0-exp}}}=|\Gal(\Qbar/E)\cdot s|^{\frac{\refC{Param0-exp}}{\refC{LGO0-exp}}}.\]
Comparing height functions and applying Pila--Wilkie (Theorem \ref{teo:PW}), we conclude that, for any $0<\epsilon<1$, there exists $\newC{PW1}:=\refC{PW1}(\cV,\epsilon)>0$ such that the number of such $x_\sigma$ lying on $\cV\setminus\cV^{\rm alg}$ is at most $\refC{PW1}|\Gal(\Qbar/E)\cdot s|^\epsilon$ (which is less than $|\Gal(\Qbar/E)\cdot s|$ when $|\Gal(\Qbar/E)\cdot s|$ is sufficiently large). Therefore, when $\Delta(s)$, and hence $|\Gal(\Qbar/E)\cdot s|$, is sufficiently large, there exists $\sigma(s)\in\Gal(\Qbar/E)\cdot s$ such that $x_\sigma\in\cV^{\rm alg}$. By \cite[Lem. B.3]{kuy:ax-lindemann}, this $x_\sigma$ is contained in a positive-dimensional irreducible algebraic subvariety of $X$ contained in $\pi^{-1}(V)$, and so, by Ax--Lindemann, $\sigma(s)$ is contained in a positive-dimensional weakly special subvariety $Z$ contained in $V$. Since $Z$ contains a special point, it is a special subvariety. By \cite[Cor. 5.3]{Orr:Galois}, we conclude that $s$ is contained in the positive-dimensional special subvariety $\sigma^{-1}(Z)$, which is also contained in $V$. Since there are only finitely many special points $s$ for which $\Delta(s)$ is less than a given bound, the proof is complete.

\begin{remark}
    The effective Chebotarev estimate supplied by Frei--Sofos--Widmer referred to in Remark \ref{KUY:unconditional} also facilitates an unconditional proof of Andr\'e--Oort via the Pila--Zannier strategy. Indeed, Ullmo--Yafaev's conditional bound \cite[Th\`eor\'eme 6.1]{uy:galois} can be made unconditional using \cite[Corollary 1.2]{FSW}. This was Conjecture 22.1 in \cite{Daw15}, which describes an otherwise unconditional implementation of the strategy, since heights of pre-special points are bounded by \cite[Th. 1.4]{DO16}.
\end{remark}

\subsection{Zilber--Pink}

The Pila--Zannier strategy seems, at present, the most promising for the purposes of attacking the Zilber--Pink conjecture. It was first described for Zilber--Pink in $Y(1)^n$ by Habegger--Pila \cite{HP16}, and in general Shimura varieties by Ren and the present author \cite{DR18}. This version of the strategy is a natural generalisation of the one just described for Andr\'e--Oort. In particular, it can also be decomposed into geometric and arithmetic components. For the former, we require the following definition. 

\begin{definition}
    Let $S$ be a Shimura variety and let $V$ be an irreducible subvariety. An irreducible subvariety $W$ of $V$ is a {\it weakly optimal subvariety} of $V$ if, for any irreducible subvariety $Y$ of $V$ strictly containing $W$, we have $\delta_{\rm ws}(Y)>\delta_{\rm ws}(W)$.
\end{definition}

The geometric component can now be stated as follows.

\begin{theorem}[Geometric Zilber--Pink, {\cite[Prop. 6.3]{DR18}}]\label{teo:GZP}
    Let $S$ be a Shimura variety and let $V$ be an irreducible subvariety of $S$. Then there exists a finite set $\Theta$ of special subvarieties of $S$ such that, for any weakly optimal subvariety $W$ of $V$, the weakly special subvariety $\langle W\rangle_{\rm ws}$ is a fibre of some $Z\in\Theta$.
\end{theorem}

Alongside o-minimality, the main ingredient in the original proof of Theorem \ref{teo:GZP} is the so-called {\it Ax--Schanuel} theorem for Shimura varieties, which was also first proved using o-minimality by Mok--Pila--Tsimerman.

\begin{theorem}[Ax--Schanuel for Shimura varieties {\cite{MPT:AS}}]
    Let $S=\Gamma\backslash X$ be a connected component of a Shimura variety and let $D$ denote the graph of the uniformisation map $\pi:X\to S$ in $X\times S$. Let $V$ be a subvariety of $X\times S$ and let $U$ be an (analytic) irreducible component of $V\cap D$. If
    \begin{equation}\label{eq:AxS}
        \dim U>\dim V-\dim S
    \end{equation}
    then the projection of $U$ to $S$ is contained in a proper weakly special subvariety.
\end{theorem}

The inequality (\ref{eq:AxS}) says that the dimension of $U$ is larger than expected. Indeed, the intersection of $D$ and $V$ in $X\times S$ would be expected to have dimension
\[\dim D+\dim V-\dim X\times S=\dim V-\dim S.\]

Ax--Schanuel was proved for $S=Y(1)^n$ by Pila--Tsimerman \cite{PT:AxS-j}. It was proved for mixed Shimura varieties of so-called {\it Kuga type} by Gao \cite{Gao:AxS}, who employed it in his aforementioned study of Betti maps. It was proved for variations of Hodge structures by Bakker--Tsimerman \cite{BT:AxS}. For variations of mixed Hodge structures, see \cite{Chiu:AxS,Chiu:AxSII} and \cite{GK:AxS}.

A refined version of geometric Zilber--Pink is obtained in \cite[Th. 3]{BD}. The method of proof is also different, working in the so-called {\it standard principal bundle} (see \cite[Sec. 5.1]{BD}) and applying multiplicity estimates to the canonical foliation thereon. This approach yields effective results (see \cite[Th. 5, 6, 7]{BD}.) The proof also employs Ax--Schanuel, which, at around the same time, was also proved in the setting of principal bundles, using differential algebra \cite{ASbundles}.  

The geometric Zilber--Pink conjecture was first proved for $Y(1)^n$ by Habegger--Pila \cite[Prop. 6.6]{HP16}. It has since been proved for variations of Hodge structures by Baldi--Klingler--Ullmo \cite[Th. 6.1]{BKU} and mixed Shimura varieties of Kuga type by Gao \cite[Th. 8.2]{Gao:AxS}. Baldi--Urbanik have recently established the conjecture for variations of mixed Hodge structures \cite[Th. 1.9]{BU:effective}. The latter proof is closer in spirit to \cite{BD} and describes an explicit algorithm.

\quad

It is relatively straightforward, using geometric Zilber--Pink, to show that the Zilber--Pink conjecture (for all $S$) is equivalent to the following (for all $S$) (see \cite[Th. 8.3]{DR18}).

\begin{conjecture}[Arithmetic Zilber--Pink]\label{conj:AZP}
    Let $S$ be a Shimura variety and let $V$ be an irreducible subvariety of $S$. Then $\Opt_0(V)$ is finite.
\end{conjecture}

When $V$ is not defined over $\Qbar$, one can sometimes exploit non-arithmetic techniques to prove cases of Conjecture \ref{conj:AZP} (see \cite[Th. 1.1]{BD:Hecke}, \cite{KT:notQbar} and \cite{AEF:Qbar}). Thus, in order to isolate the arithmetic aspects, we will henceforth assume that $V$ is defined over $\Qbar$.    

In this case, as for Andr\'e--Oort, the Pila--Zannier strategy for Conjecture \ref{conj:AZP} requires two additional ingredients. In order to state them, we revive the notation from Section \ref{subsec:AO} and let $\Sigma$ denote the set of special subvarieties of $S$. Extend $\Delta$ to a complexity $\Sigma\to\bN$ again with the property that, for any $B>0$, the set
\[\{Z\in\Sigma:\Delta(Z)<B\}\]
is finite. Furthermore, for simplicity, choose $\Delta$ such that every special subvariety $Z$ contains a special point $s$ satisfying $\Delta(s)\leq\Delta(Z)$ (cf. \cite[Def. 10.2]{DR18}). For brevity, we now say {\it definable} instead of {\it definable in $\bR_{\rm an,exp}$}.

\subsubsection{Large Galois orbits}\label{subsec:LGO}
The first of the ingredients is to show that there exists $\newC{LGO-exp}>0$ such that, for any $s\in\Opt_0(V)$, we have
\[|\Gal(\Qbar/E)\cdot s|\gg\Delta(\langle s\rangle)^{\refC{LGO-exp}}.\]
This problem is still wide open, and, seemingly, very difficult. The cases hitherto known (for $s$ non-special) are those described in Section \ref{sec:ZPhistory} (though the restriction to {\it simple} PEL type appearing there is due to the state of the art towards the next ingredient rather than this one). Except for the result for asymmetric curves, these cases have been obtained using techniques originally due to Andr\'e \cite{And89} and Bombieri \cite{bombieri:Gfns}, based on the theory of G-functions. This approach will be described in Section \ref{sec:Gfns}.

\begin{remark}
    Lower bounds for Galois orbits of points in Hecke orbits are the main technical achievement of \cite{RY:APZ}, which also uses the Pila--Zannier strategy. For an outline of the strategy in that setting, see \cite[Sec. 1.4]{RY:heights}.
\end{remark}

\subsubsection{Parameter height bounds}\label{subsec:PHBs}
The second ingredient is to construct a definable parameter space $\cM\subset\bR^d$ for totally geodesic irreducible algebraic subvarieties of $X$ and to find $k\in\bN$ and $\newC{PHBs-exp}>0$ such that, for $s\in S$, there exists $m\in\cM(\Qbar):=\cM\cap\Qbar^d$ representing the totally geodesic subvariety $Y_m$ such that
    \begin{enumerate}[(i)]
        \item $Y_m\cap\pi^{-1}(s)\cap\cF\neq\emptyset$;
        \item $\pi(Y_m)=\langle s\rangle$;
        \item ${H}_k(m)\ll\Delta(\langle s\rangle)^{\refC{PHBs-exp}}$.
    \end{enumerate}

Progress towards this ingredient is more advanced. The result for special points (in which case, $\cM=X$) is \cite{DO16} (see Section \ref{sec:hts-pre-sp}). The result for special subvarieties of $Y(1)^n$ follows immediately from \cite[Lem. 5.2]{HP12}. Orr generalised this result to the Hecke orbit of a special subvariety in an arbitrary Shimura variety \cite[Lem. 3.9]{OrrHecke} (based on the main result of \cite{Orr18}), and this was generalised to so-called Hecke-factor families in \cite[Lem. 6.1]{ExCM}. In \cite{QRTUI}, a general approach towards parameter height bounds was introduced (see Section \ref{sec:PHBs}) and applied to so-called quaternionic and $E^2$ curves in $\cA_2$ \cite[Prop. 6.3]{QRTUI}. The paper \cite{PEL} handles special subvarieties of $\cA_g$ parametrising abelian varieties whose endomorphism ring contains a simple division algebra of type I or II in Albert's classification (see \cite[\S21, Th. 2]{mumford:abelian}). Types III and IV are treated in the work of Bhatta \cite{Bhatta}. A new approach due to Urbanik (see \cite[Prop. 3.2]{urbanik:CASP}, for example) applies to general variations of Hodge structures.

\subsubsection{Pila--Zannier for Arithmetic Zilber--Pink}

Suppose $s\in\Opt_0(V)$. This implies that $s$ is a component of $V\cap\langle s\rangle$. By assumption, $V$ is defined over $E$, and so $s$ is defined over $\Qbar$. By \cite[Cor. 5.3]{Orr:Galois}, for every $\sigma\in\Gal(\Qbar/E)$, we have $\sigma(s)\in\Opt_0(V)$. 

For each $\sigma(s)\in\Gal(\Qbar/E)\cdot s$, we let $m_\sigma\in\cM(\Qbar)$ be the point afforded to us by Section \ref{subsec:PHBs} and write $Y_\sigma:=Y_{m_\sigma}$ for the corresponding totally geodesic subvariety. In particular, we have $\pi(Y_\sigma)=\langle \sigma(s)\rangle$ and, writing $\Delta_\sigma:=\Delta(\langle\sigma(s)\rangle)$, 
\[{H}_k(m_\sigma)\ll\Delta^{\refC{PHBs-exp}}_\sigma\ll|\Gal(\Qbar/E)\cdot \sigma(s)|^{\frac{\refC{PHBs-exp}}{\refC{LGO-exp}}}=|\Gal(\Qbar/E)\cdot s|^{\frac{\refC{PHBs-exp}}{\refC{LGO-exp}}},\]
by Sections \ref{subsec:LGO} and \ref{subsec:PHBs}.
Moreover, we may choose $z_\sigma\in Y_\sigma\cap\pi^{-1}(\sigma(s))\cap \cF$.

The points $(m_\sigma,z_\sigma)$ belong to the definable set
\[\cD:=\{(m,z)\in\cM\times\cV:z\in Y_m\}.\]
Therefore, by the semi-rational Pila--Wilkie theorem (Theorem \ref{teo:PWSR}) applied to $\cD$, if $\Delta(\langle s\rangle)$ is sufficiently large, there exists a continuous definable function $\beta:[0,1]\to\cD$, which is real analytic on $(0,1)$, whose projection to $\cM$ (resp. $\cV$) is semi-algebraic (resp. non-constant), and satisfies $\beta(0)=(m_\sigma,z_\sigma)$ for some $\sigma\in\Gal(\Qbar/E)$ (for the last assertion, we appeal to the slight strengthening \cite[Th. 9.1]{DR18}). Using o-minimality and Ax--Schanuel, Cassani has shown in his PhD thesis \cite{cassani:thesis} that this leads to a contradiction. (For the case $\dim V=1$, see \cite[Th. 14.2]{DR18}.) 

We conclude that $\Delta(\langle s\rangle)$ is bounded and, therefore, $\langle s\rangle$ belongs to a finite set. This concludes the argument (again using that $s$ is an irreducible component of $V\cap\langle s\rangle$).

\section{Large Galois orbits: G-functions}\label{sec:Gfns}
 
In this section, we will loosely summarise the aforementioned approach towards large Galois orbits for the case $S=\cA_g$ and unlikely intersections with PEL type special subvarieties.   

At present, the method can only handle cases for which $\dim V=1$. Therefore, we rewrite $C:=V$. After some minor adjustments, $C$ is the base of an abelian scheme $f:\cA\to C$ of relative dimension $g$ defined over a number field $K\subset\bC$. To this, we can associate the first relative algebraic de Rham cohomology $H^1_{\rm dR}(\cA/C)$, which is equipped with an integrable $\cO_C$-connection
\[\nabla:H^1_{\rm dR}(\cA/C)\to H^1_{\rm dR}(\cA/C)\otimes_{\cO_C}\Omega^1_C,\]
known as the {\it Gauss--Manin connection} on $\cA\to C$ \cite{KO:GM}. Since $H^1_{\rm dR}(\cA/C)$ is locally free, we can, after deleting finitely many points from $C$, pick an isomorphism $\cO^{2g}_C\cong H^1_{\rm dR}(\cA/C)$ and let $\{\omega_i\}_{i=1}^{2g}$ denote the image of the standard (global) basis. On the other hand, for a sufficiently small open set $\Delta\subset C^{\rm an}$ we can choose a basis $\{\gamma_i\}_{i=1}^{2g}$ for $H^1_{\rm dR}(\cA/C)^{\rm an}(\Delta)=H^1_{\rm dR}(\cA^{\rm an}/C^{\rm an})(\Delta)$ that is horizontal for $\nabla^{\rm an}$ and let $\Omega\in\rM_{2g}(\cO_{C^{\rm an}}(\Delta))$ denote the matrix whose rows give the coordinates of $\{\omega^{\rm an}_i\}_{i=1}^{2g}$ in terms of $\{\gamma_i\}_{i=1}^{2g}$. Note that, if we imposed that the $\gamma_i$ belong to the image of $R^1\!f^{\rm an}_*\,\Qbar$ under the map \cite[Ch. IX, Sec. 1.2 (1.2.2)]{And89}, it would be usual to refer to the entries of $\Omega$ as {\it periods} of $\cA\to C$.  

We will base our constructions around a fixed point $s_0\in C(K)$ and work with a local parameter $x\in K(C)$ at $s_0$. In most of the work to date, it is actually necessary to enlarge $C$ and choose $s_0$ to be a point for which $\cA\to C$ has multiplicative degeneration. It is then necessary to extend various objects across $s_0$. For simplicity, however, we will suppress this complication. Our work with Orr and Papas \cite{DOP}, for example, handles some cases for which $s_0$ is a genuine point of $C$.  

We define the normalised matrix $Y:=\Omega\cdot\Omega(s^{\rm an}_0)^{-1}$ and define $(F_{ij})_{i,j=1}^{2g}$ to be the matrix whose entries are the Taylor series of the entries of $Y$ with respect to $x$. It is possible to show that the $F_{ij}$ are a particular type of power series, known as {\it G-functions}, introduced by Siegel in his 1929 paper (see \cite[p. 1]{And89} for the definition to which we refer). In particular, they satisfy a linear homogeneous differential equation induced by $\nabla$, and they  belong to $K[[X]]$, after possibly replacing $K$ by a finite extension.

The Hasse principle of Bombieri \cite[Intro., Th. E.]{And89} (cf. \cite[Sec. 11]{bombieri:Gfns}) states that, for a collection of G-functions $F_1,\ldots,F_n$, there exists $\newC{mult-Gfn}>0$ such that the elements of
\[\{x\in\Qbar:F_1,\ldots,F_n\text{ satisfy a non-trivial global relation of degree }d\text{ at }x\}\]
satisfy $h(x)\ll d^{\refC{mult-Gfn}}$. A relation refers to a homogeneous polynomial over $\Qbar$ in $n$ variables. We say that $F_1,\ldots,F_n$ satisfy a relation $P$ at $x$ for a place $v$ of $K$ (after possibly replacing $K$ with a finite extension containing $x$ and the coefficients of $P$) if $|x|_v<1$, the $F_i$ converge $v$-adically at $x$, and $P(F_1(x),\ldots,F_n(x))=0$ in $K_v$.  We say that $P$ is a global relation if $F_1,\ldots,F_n$ satisfy $P$ at $x$ for all places $v$ for which $|x|_v<1$ and the $F_i$ converge $v$-adically at $x$. A global relation is non-trivial if it is not the specialisation of a functional relation between the $F_i$ (see \cite[Ch. 7, Sec. 5.1]{And89}). 

The strategy towards large Galois orbits is to show that, at a point $s\in C(\Qbar)$ belonging to a PEL type special subvariety of $\cA_g$, there exists a global non-trivial relation between the $F_{ij}$ at $x(s)$ whose degree is bounded by a positive power of $[K(s):K]$. In this case, Bombieri yields
\[h(s)\ll [K(s):K]^{\newC{Andre-mult}}\]
and this bound can be combined with the results of Masser--W\"ustholz \cite{MW93,MW:endo} to yield lower bounds for Galois orbits (see, for example, the proof of \cite[Th. 6.5]{QRTUI}).

This technique first appeared in \cite{And89} and was further utilised in \cite{And95}. As alluded to above, it was revived in \cite{ExCM} and has since been extended in the abelian setting in the aforementioned papers of Orr, Papas and the present author (see Section \ref{sec:ZPhistory}) and in \cite{Youell}. As mentioned in Section \ref{sec:AOhistory}, Binyamini--Masser independently revived Andr\'e's techniques in their aforementioned work on effective Andr\'e--Oort in Hilbert modular varieties \cite{BM:AO}. In \cite{Urbanik:Gfns}, Urbanik generalises the approach to variations of Hodge structures. For a wider survey of G-functions in arithmetic geometry, see the recent survey of Andr\'e \cite{Andre:survey}.

\section{Parameter height bounds: quantitative reduction theory}\label{sec:PHBs}

In this section, we describe an existing approach towards parameter height bounds. To that end, let $S$ be a Shimura variety defined by $(\gG,\gX)$ and let $\Gamma\backslash X$ denote one of its connected components.

By definition, special subvarieties of $\Gamma\backslash X$ are the images of those subvarieties $X_\gH$ of $X$ with $(\gH,\gX_\gH)$ a Shimura subdatum of $(\gG,\gX)$ and $X_\gH$ a connected component of $\gX_\gH$ contained in $X$. For such an $X_{\gH}$, we have
\[X_\gH=\gH(\bR)^+x=\gH^{\rm der}(\bR)^+x\]
for any $x\in X_\gH$ \cite[Prop. 5.7(a)]{Mil05}, and so we can parametrise special subvarieties (non-uniquely) using pairs $(\gH,x)$ with $\gH$ a semisimple $\bQ$-subgroup of $\gG$ and $x\in X$ such that $x$ factors through the almost direct product $\gH_\bR\gZ_\gG(\gH)_\bR$. (Note that, for a reductive subgroup $\gH$ of $\gG$, we have $\gH=\gH^\der\gZ(\gH)\subset\gH^\der\gZ_\gG(\gH^\der)$.) These pairs form a subset of the pairs $(\gH,x)$ with $\gH$ a semisimple $\bR$-subgroup of $\gG_\bR$ possessing no compact factors and $x\in X$ such that $x$ factors through $\gH\gZ_{\gG_\bR}(\gH)$. The latter correspond exactly to the totally geodesic subvarieties of $X$ \cite[Prop. 2.3]{uy:algebraic-flows}.

To parametrise a special subvariety $Z$, we choose $x\in\cF$ pre-special such that $\Delta(\pi(x))\leq\Delta(Z)$. Thus, it remains to parametrise the semisimple $\bQ$-subgroups of $\gG$ associated with special subvarieties of $S$. Since these belong to finitely many $\gG(\bR)$-conjugacy classes (see \cite[Cor. 0.2]{BDR}, for example), it suffices to work in the conjugacy class of a given semisimple $\bQ$-subgroup $\gH_0=\gMT(\gX_0)^\der$ associated with a Shimura subdatum $(\gMT(\gX_0),\gX_0)$ of $(\gG,\gX)$. To simplify the exposition (cf. \cite[Cor. 5.3]{PEL}), we will assume that there exists a connected component $X_0$ of $\gX_0$ such that, for any subvariety $X_\gH$ as above, with $\gH^\der_\bR=g\gH_{0,\bR}g^{-1}$ for some $g\in\gG(\bR)$, we have \[X_\gH=gX_0=g\gH_0(\bR)^+x_0=g\gH_0(\bR)^+g^{-1}\cdot gx_0.\]

One approach is as follows. It is straightforward to construct a finite-dimensional $\bQ$-representation $\gG\to\gGL(V)$ and a vector $v_0\in V$ such that $\gH_0$ is the stabiliser $\Stab_\gG(v_0)$ of $v_0$ in $\gG$ (see \cite[Th. 16.1]{Waterhouse}). And, of course, we may choose a free $\bZ$-module $\Lambda\cong\bZ^n$ such that $V=\Lambda_\bQ$ and $v_0\in\Lambda$. In \cite{QRTUI}, we show that, for such a representation (provided $\gG(\bR)v_0$ is closed in $V_\bR$), there exists $\newC{PHBs-exp1}>0$ such that, for any $g\in\gG(\bR)$ and $v_g\in\gZ_{\gGL(V)}(\gG)(\bR)v_0$ satisfying
\begin{itemize}
\item[(1)] $\gH_g:=g\gH_{0,\bR}g^{-1}$ is defined over $\bQ$ and 
\item[(2)] $gv_g\in\Lambda$,
\end{itemize}
there exists a fundamental set for $\Gamma_g:=\Gamma\cap\gH_g(\bR)$ in $\gH_g(\bR)$ of the form
\[B_g\cF_{\gG} g^{-1}\cap\gH_g(\bR)\]
for some finite $B_g\subset\Gamma$ with the property that $|b^{-1}gv_g|\ll|v_g|^{\refC{PHBs-exp1}}$ for every $b\in B_g$. It is this aspect of the strategy from which the name {\it quantitative reduction theory} comes.

Now we imitate the proof of \cite[Lem. 8.3]{PEL}. For $s\in S$ and $Z_g:=\pi(gX_0)=\langle s\rangle$, we pick $y\in gX_0\cap\pi^{-1}(s)$. By the above, we can write $y=\gamma bfg^{-1}\cdot gx_0$ for some $\gamma\in\Gamma_g$, $b\in B_g$ and $f\in\cF_{\gG}$, and we set
\[z=b^{-1}\gamma^{-1}y=fx_0\in\cF_{\gG} x_0\cap X=\cF.\]
Since $b,\gamma\in\Gamma$, we have $z\in\pi^{-1}(s)$. Moreover, $z$ belongs to the totally geodesic subvariety $b^{-1}gX_0$ parametrised by $b^{-1}gv_g\in\Lambda\subset V_\bR$ (whose image in $S$ is $Z_g$).

Therefore, in order to obtain parameter height bounds via this approach, it remains to construct representations and $v_g$ as above satisfying $|v_g|\ll\Delta(Z_g)^{\newC{PHBs2-exp}}$ for some $\refC{PHBs2-exp}>0$ independent of $g$. This is carried out in \cite{QRTUI}, \cite{PEL} and \cite{Bhatta} for simple PEL type special subvarieties.

\bibliographystyle{plain}
\bibliography{ZP}

\end{document}